\documentclass[12pt,a4paper]{article}

\usepackage{amssymb}
\usepackage{yfonts}
\usepackage{amscd,diagrams}

\usepackage[T1]{fontenc}
\def\og{\guillemotleft\,}
\def\fg{\,\guillemotright}

\DeclareMathSymbol{\Z}{\mathalpha}{AMSb}{"5A} 
\DeclareMathSymbol{\PP}{\mathalpha}{AMSb}{"50} 
\DeclareMathSymbol{\Q}{\mathalpha}{AMSb}{"51}
\DeclareMathSymbol{\N}{\mathalpha}{AMSb}{"4E}
\DeclareMathSymbol{\R}{\mathalpha}{AMSb}{"52}

\newtheorem{theo}{Th{\'e}or{\`e}me}
\newtheorem{Theo}{Th{\'e}or{\`e}me}

\newtheorem{prop}{Proposition}
\newtheorem{lem}{Lemme}
\newtheorem{cor}{Corollaire}
\newenvironment{rem}{\noindent{\bf Remarque.}}{}

\newenvironment{dem}{\noindent {\sc D{\'e}monstration.}}{\hfill$\square$}
\newenvironment{preuve}{\noindent {\sc Preuve.}}{\hfill$\square$}

\newcounter{exple} 
\newenvironment{exple}{\refstepcounter{exple}\noindent{\bf Exemple \arabic{exple}.}
}{}

\newcommand\Ker[2]{{\rm Ker}\left ( #1\longto{\rm Aut}(#2)\right )}

\def\ouv{\Omega}

\def\modeleB1{\hpu_1/\gp^u_1\times \hD_1/\hBD}
\def\deltaF{{\delta_{\mathcal F}}}
\def\cF{{\mathcal F}}
\def\gp{{\textgoth{p}}}
\def\longto{\longrightarrow}
\def\hG{\hat{G}}
\def\hP{\hat{P}}
\def\hD{\hat{D}}

\def\hB{\hat{B}}

\def\hpu{\hat{\gp_u}}

\def\hBD{\hat{B_D}}
\def\DQG{\backslash\hspace*{-1.5pt}\backslash}

\def\DQD{/\hspace*{-1.5pt}/}

\def\CF{\mathcal C_{\mathcal F}}
\def\CFun{\mathcal C_{\mathcal F_1}}
\def\CFdeux{\mathcal C_{\mathcal F_2}}

\newcommand{\mc}[1]{\mathcal{#1}}

\newcommand{\extdeux}{\Lambda^2}

\newcommand{\h}[1]{\hat{#1}} 
\newcommand{\psd}{\rtimes}

\newcommand{\End}{\mathrm{End}}
\newcommand{\Hom}{\mathrm{Hom}}

\begin{document}

\title{Sur des faces du LR-c{\^o}ne g{\'e}n{\'e}ralis{\'e}}
\author{Pierre-Louis Montagard et Nicolas Ressayre}

\maketitle

\section{R{\'e}sumé}
Soient $G\subset\h G$ deux groupes r{\'e}ductifs connexes. 
Notons $\mc D$ (resp. $\h{\mc D}$)  l'ensemble des classes d'isomorphisme des
repr{\'e}sentations irr{\'e}ductibles de $G$ (resp. $\h G$). 
Nous nous int{\'e}ressons {\`a} l'ensemble
$\mathcal{C}$ des couples $(\mu,\h\nu)$ dans $\mc D\times\h{\mc D}$ 
pour lesquels un $\h G$-module
de classe $\h \nu$ contient un sous-$G$-module de classe $\mu$.
Il est bien connu que $\mathcal{C}$ engendre un c{\^o}ne poly{\'e}dral
dans un espace vectoriel appropri{\'e}.
Par des m{\'e}thodes de th{\'e}orie g{\'e}om{\'e}trique des 
invariants nous {\'e}tudions sous quelles conditions une in{\'e}galit{\'e} 
lin{\'e}aire
d{\'e}finissant $\mc D$ induit une face de codimension un du c{\^o}ne
engendr{\'e} par $\mathcal{C}$.

\section{Introduction}\label{intro}

Soit $k$ un corps alg{\'e}briquement clos de caract{\'e}ristique nulle.
Soient $G\subset \h{G}$ deux groupes alg{\'e}briques
r{\'e}ductifs connexes sur $k$.
Soit $\h{V}$ une repr{\'e}sentation rationnelle de
dimension finie de $\h G$. Le probl{\`e}me g{\'e}n{\'e}ral que l'on aborde ici est
de d{\'e}composer $\h V$ en somme de $G$-modules irr{\'e}ductibles.
Remarquons que ce contexte recouvre de nombreux probl{\`e}mes de
d{\'e}composition de 
repr{\'e}sentations. 
Citons en deux~:
\begin{itemize}
\item\label{produittens}
si $\h G=G\times G$ et
si $G$ est la diagonale de $\h G$, il s'agit de d{\'e}composer le
produit tensoriel de deux repr{\'e}sentations irr{\'e}ductibles de $G$ ;

\item\label{plethys}
soit $G$ un groupe r{\'e}ductif quelconque et 
$\rho\ :\ G\rightarrow {\rm Gl}(V)$ une repr{\'e}sentation
irr{\'e}ductible de $G$, on peut alors poser $\h G:={\rm Gl}(V)$ et
consid{\'e}rer l'inclusion $\rho(G)\subset \h G$, le probl{\`e}me est
alors de d{\'e}composer des repr{\'e}sentations de $G$ telles
que  la puissance sym{\'e}trique $n$-i{\`e}me $S^nV$ de $V$,
 la puissance ext{\'e}rieure $k$-i{\`e}me $\Lambda^kV$ de $V$, ou plus
g{\'e}n{\'e}ralement de d{\'e}composer les puissances de Schur $S_\pi V$.
\end{itemize}

Dans ce contexte tr{\`e}s g{\'e}n{\'e}ral, on ne cherche pas {\`a} donner des
formules combinatoires explicites de d{\'e}composition comme la
c{\'e}l{\`e}bre r{\`e}gle de Littlewood-Richardson concernant la
d{\'e}composition du produit tensoriel pour le groupe lin{\'e}aire. Notre
approche est plus qualitative {\`a} travers le c{\^o}ne de
Littlewood-Richardson g{\'e}n{\'e}ralis{\'e} que nous allons d{\'e}finir,
apr{\`e}s avoir introduit quelques notations suppl{\'e}mentaires.

Nous noterons $\mathcal{D}$ (resp. $\h{\mathcal{D}}$)
l'ensemble des classes d'isomorphismes de repr{\'e}sentations
irr{\'e}ductibles de $G$ (resp. $\h G$). Pour $\nu\in\mc D$
(resp. $\h\nu\in \h{\mc D}$) nous noterons $V_\nu$ (resp. $V_{\h\nu}$)
une repr{\'e}sentation irr{\'e}ductible de $G$ (resp. $\hG$) dans la
classe $\nu$ (resp.  $\h\nu$). 
Si $V$ est une repr{\'e}sentation de $G$ nous noterons $(V_\nu,V)$ la 
multiplicit{\'e} de
$V_\nu$ dans $V$, c'est-{\`a}-dire la dimension de l'espace des homomorphismes
$G$-{\'e}quivariants de $V_\nu$ vers $V$. Nous nous int{\'e}ressons
{\`a} 
l'ensemble~:

$$\mathcal{C}:=\{(\mu,\h\nu)\in\mathcal{D}\times\h
{\mathcal{D}}\,|\,(V_\mu, V_{\h \nu})\neq 0\}\,.$$
 
Les ensembles $\mc D$ et $\hat\mc D$ ont une structure naturelle de
semi-groupes. De plus, M.~Brion et F.~Knop ont montr{\'e} (voir
\cite{elash}) que $\mc C$ est un sous-semigroupe de type fini de
$\mathcal{D}\times\hat\mathcal{D}$. Dans le cas du produit tensoriel
pour le groupe lin{\'e}aire ce semi-groupe  a {\'e}t{\'e} appel{\'e}
semi-groupe de Littlewood-Richardson (voir \cite{zelevinski}).
Comme $\mathcal{D}$ (resp. $\h{\mathcal{D}}$) est en bijection avec les
points entiers d'un c{\^o}ne d'un $\Q$-espace vectoriel que nous
appellerons provisoirement $E$ (resp. $\h E$),
$\mc C$ engendre un c{\^o}ne poly{\'e}dral $\tilde{\mathcal C}$ dans 
$E\times \h E$. Nous appelons ce dernier 
{\it c{\^o}ne de Littlewood-Richardson g{\'e}n{\'e}ralis{\'e}}, ou plus 
bri{\`e}vement 
{\it LR-c{\^o}ne g{\'e}n{\'e}ralis{\'e}}.  
Comme  $\tilde{\mathcal C}$ est poly{\'e}dral, il est d{\'e}fini dans 
$E\times \h E$ par un nombre fini d'in{\'e}galit{\'e}s lin{\'e}aires 
correspondantes aux faces de codimension un, 
faces que nous appellerons essentielles. 
Par la suite, pour une partie $\mc E$ d'un espace vectoriel 
nous appellerons dimension de $\mc{E}$ et noterons dim$\mc{E}$  
la dimension de l'espace vectoriel engendr{\'e} par $\mc E$.
Soit $\mathcal{F}$ l'espace vectoriel engendr{\'e} par une face du
c{\^o}ne engendr{\'e} par $\mc D$. Il induit naturellement une \og face\fg\ 
$\mc C_{\mc F}$ de $\mc C$ d{\'e}finie par~:
$\mc C_{\mc F}:=\mc C\cap(\mc F\times\h E)$. Le but de cet article est
d'{\'e}tudier la dimension de $\mc C_{\mc F}$ et notamment de savoir 
si elle engendre une face essentielle de $\tilde{\mathcal C}$. 
Remarquons que dans le cas du produit tensoriel
pour le groupe lin{\'e}aire $\mathrm{GL}(V)$, Knutson, Terao et Woodward
\cite{KTW} d{\'e}crivent toutes les faces essentielles de $\tilde{\mc C}$, 
et montrent notamment que si $\mc F$ est une face essentielle de $\mc D$,
alors elle induit une face essentielle de $\tilde{\mc C}$ d{\`e}s que la dimension
de $V$ est sup{\'e}rieure ou {\'e}gale {\`a} $3$.
Nous retrouvons et {\'e}largissons ce r{\'e}sultat dans la 
proposition~\ref{exple:GGG}.

Dans le cas g{\'e}n{\'e}ral, des
consid{\'e}rations {\'e}l{\'e}mentaires (voir le lemme~\ref{dimClin}) permettent
de montrer que la codimension de
$\CF$ dans $\mc C$ est sup{\'e}rieure ou {\'e}gale {\`a} la
codimension de $\mc F$ dans $\mc D$. 
Cette in{\'e}galit{\'e} se traduit par 
$\deltaF:=\dim\mc F-\dim\CF-\dim\mc D+\dim\mc C\geq 0$. 
On dit que $\mc F$ est {\it pleine} si $\deltaF=0$, c'est-{\`a}-dire
si dim$\mc C_{\mc F}$ est maximale.
En particulier, si $\mc F$ provient d'une face essentielle de 
$\mc D$, alors $\deltaF=0$ est {\'e}quivalent au fait que  
$\CF$ soit essentielle dans $\mc C$. 
On peut maintenant {\'e}noncer un de nos r{\'e}sultats~:
(voir le corollaire~\ref{cor:monotonie})

\begin{Theo}\label{mainthm1} 
Soit $\mc F_1$ et $\mc F_2$ deux sous-espaces de $E$ engendr{\'e}s par deux 
faces du c{\^o}ne engendr{\'e} par $\mc D$.
Si $\mc F_1\subset \mc F_2$, alors $\deltaF_1\geq \deltaF_2$.
\end{Theo}

Nous pr{\'e}sentons {\'e}galement dans ce travail une condition
{\'e}quivalente au fait que $\mc F$
soit pleine. Pour pouvoir {\'e}noncer celle-ci, nous allons introduire
quelques d{\'e}finitions suppl{\'e}mentaires.

De mani{\`e}re classique (voir la section~\ref{sec:noyau}),  on associe 
{\`a} $\mc F$ un groupe parabolique $P$ de $G$ ; le groupe $P$ se d{\'e}compose en un produit
semi-direct $P^u\psd L$, o{\`u} $P^u$ est le radical unipotent de $P$
et $L$ un sous-groupe de L{\'e}vi. 
Dans la section~\ref{choix}, nous montrons alors qu'il existe un
sous-groupe parabolique $\h P$ de $\h G$ et une
d{\'e}composition de L{\'e}vi de celui-ci~: $\h P=\h P^u\psd\h L$ qui
v{\'e}rifient~: $P=\h P\cap G$ ; $P^u=\h P^u\cap G$ et $L=\h L\cap G$.

Nous noterons $\gp^u$ (resp. $\h\gp^u$),
l'alg{\`e}bre de Lie de $P^u$ (resp. $\h P^u$) et $B_L$ (resp. $B_{\h L}$) un
sous-groupe de Borel de $L$ (resp.  $\h L$). 
Enfin rappelons que si un groupe alg{\'e}brique $\Gamma$ agit sur une
vari{\'e}t{\'e} $X$, on appelle isotropie r{\'e}ductive de $\Gamma$ en
$x\in X$, le quotient de $\Gamma_x$ par son radical unipotent.
Le groupe $\h L$ agit sur $\h\gp^u$ par la repr{\'e}sentation adjointe 
et $L\subset\h L$ stabilise $p^u$~; 
donc, $L$ agit sur $\h\gp^u/\gp^u$. 
De plus, $\h L$ et donc $L$ agissent sur $\h L/B_{\h L}$ par multiplication. 
Finalement, $L$ agit sur $\h\gp^u/\gp^u\times\h L/B_{\h  L}$ diagonalement.
On peut maintenant {\'e}noncer le (voir le corollaire~\ref{cor:ppal})~:

\begin{Theo}{\label{mainthm2}}
Il existe un ouvert non vide $\ouv$ de $\h\gp^u/\gp
^u\times
\h L/B_{\h L}$ tel que pour tout $x\in\ouv$, $\deltaF$ est
{\'e}gal {\`a} la diff{\'e}rence des dimensions des isotropies r{\'e}ductives
des groupes $L$ et $B_L$ en $x$.
\end{Theo}

Une cons{\'e}quence imm{\'e}diate des th{\'e}or{\`e}mes~\ref{mainthm1} 
et~\ref{mainthm2} est le~:

\begin{Theo}\label{mainthm3}
S'il existe un point de $\h G/\h B$ dont l'isotropie
dans le groupe d{\'e}riv{\'e} de $G$ 
est finie, alors toutes les faces de $\mc D$ sont pleines.
\end{Theo}

Dans la section~\ref{sec:exple}, nous appliquons ce dernier r{\'e}sultat 
{\`a} divers exemples.

Remarquons que $\mc C$ est tr{\`e}s li{\'e} au polytope moment d{\'e}fini dans un cadre 
symplectique. La propri{\'e}t{\'e} d'\^etre pleine s'interpr{\`e}te en terme de ces 
polytopes, voir la proposition \ref{prop:polmom}.

\section{Premi{\`e}res propri{\'e}t{\'e}s}\label{premprop}
\subsection{Notations}\label{notations}

Commen{\c c}ons par quelques notations g{\'e}n{\'e}rales~:
si $\Gamma$ est un groupe alg{\'e}brique affine sur $k$, nous noterons
$\Gamma^u$ son radical unipotent, $[\Gamma,\Gamma]$ son groupe
d{\'e}riv{\'e}, $\Gamma^\circ$ sa composante neutre, 
$\Xi(\Gamma)={\rm Hom}(\Gamma,k^*)$ le groupe de ses caract{\`e}res, 
$\Xi_*(\Gamma)={\rm Hom}(k^*,\Gamma)$ le groupe de ses sous-groupes {\`a} un param{\`e}tre 
et $\mathrm{Lie}(\Gamma)$ son alg{\`e}bre de Lie. 
Dans tout cet article, nous appelons vari{\'e}t{\'e}, une vari{\'e}t{\'e} alg{\'e}brique
quasi-projective et irr{\'e}ductible. 
Si $\Gamma$ op{\`e}re alg{\'e}briquement sur une vari{\'e}t{\'e} $X$ 
on dit que $X$ est une {\it $\Gamma$-vari{\'e}t{\'e}}. 
On note $\Ker{\Gamma}{X}$ le noyau de l'action de $\Gamma$ sur $X$.
Si $X$ est affine et l'alg{\`e}bre $k[X]^\Gamma$ des fonctions r{\'e}guli{\`e}res sur $X$ 
invariantes par $\Gamma$ est de type fini, $X\DQD\Gamma$ d{\'e}signera la
vari{\'e}t{\'e} affine associ{\'e}e {\`a} $k[X]^\Gamma$. 
L'inclusion de $k[X]^\Gamma$ dans $k[X]$ induit un morphisme dominant et
$\Gamma$-invariant $\pi\,:\,X\longto X\DQD \Gamma$ que nous appelons 
{\it morphisme quotient}.
Si $V$ est un $\Q$-espace vectoriel et si $\mathcal E$ 
est un sous-ensemble de $V$,
nous noterons $<\mathcal E>$ l'espace vectoriel 
engendr{\'e} par
$\mathcal E$, $\dim\mc E$ la dimension de $<\mc E>$ et 
$\mc E^\perp$ l'ensemble des $\varphi\in V^*$ tels que $\varphi_{|\mc E}=0$. 
Enfin si $R$ est un groupe commutatif $R_\Q$ d{\'e}signera
le $\Q$-espace vectoriel $R\otimes_\Z\Q$. 

Rappelons que nous consid{\'e}rons deux groupes alg{\'e}briques 
r{\'e}ductifs connexes $G\subset\h G$.
Fixons, pour tout l'article, 
un tore maximal $T$ de $G$ et un sous-groupe de Borel $B$ de $G$.
L'ensemble $\mc{D}$ s'identifie alors naturellement au
sous-ensemble des poids dominants de $\Xi(T)$. Si $\mu\in\mc{D}$, nous
noterons $V_\mu$ une repr{\'e}sentation irr{\'e}ductible de poids
dominant $\mu$.

Soit $\h T\subset \h B$ un tore maximal et un sous-groupe de Borel de 
$\h G$. La notation $V_{\h \nu}$ d{\'e}signe une repr{\'e}sentation irr{\'e}ductible
de $\h G$ de plus haut poids $\h \nu$.

\subsection{{\'E}nonc{\'e} du probl{\`e}me}\label{enonce}

Rappelons que nous nous int{\'e}ressons {\`a} 

$$\mathcal{C}:=\{(\mu,\nu)\in\mathcal{D}\times\h
{\mathcal{D}}\,|\,(V_\mu, V_{\h \nu})\neq 0\}\,.$$

Soit $\mathcal{F}$ le sous-espace vectoriel de $\Xi(T)_\Q$ 
engendr{\'e} par une face du c{\^o}ne engendr{\'e} par $\mc D$. Dans
la suite, par abus de notation, nous appellerons $\mc F$ une face de
$\mc D$. On pose alors~: 

$$\mathcal{C}_{\mathcal{F}}:=\mathcal{C}\cap\left(\mathcal{F}\times\Xi_\Q(\hat
{T}) \right)\\
\,.
$$

Dans cet article, nous  cherchons {\`a} comparer la codimension de
$\mc F$ dans $\Xi(T)_\Q$ {\`a} celle de $\mc C_{\mc F}$ dans $\mc C$.
Pour cela on pose
$\delta_{\mc F}=\dim\mc C-\dim \mc C_{\mc F}+\dim\mc F-\dim\Xi(T)$.

\begin{lem}\label{dimClin}
On a l'{\'e}galit{\'e}~:
$$
\dim(<\mc C>\cap(\mc F\times\Xi(\h T)_{\Q}))=\dim(\mc C)-\dim(\Xi(T))+\dim(\mc F).
$$
En particulier, $\delta_{\mc F}\geq 0$.
\end{lem}

\begin{preuve} 
Posons $d=\dim \mc C -\dim{\left ((\mc F\times\Xi(\h
    T)_\Q)\cap(<\mc C>)\right )}$.
Consid{\'e}rons l'application $\pi\,:\, <\mc C>\rightarrow \Xi(T)_\Q$ induite par la
projection de $\Xi(T)\times\h\Xi(T)$ sur $\Xi(T)$ et sa transpos{\'e} 
$^t\pi\,:\,{\rm Hom}(\Xi(T)_\Q,\Q)\rightarrow
<\mc C>^*:={\rm Hom}(<\mc C>,\Q)$.

L'entier $d$ est la dimension de
l'orthogonal dans $<\mc C>^*$ de $\pi^{-1}(\mc F)$ ; ainsi,
$d= \dim(^t\pi(\mc F^\perp))$.
Or, comme toute
repr{\'e}sentation irr{\'e}ductible de $G$
appara{\^\i}t dans au moins une repr{\'e}sentation de $\h G$, $\pi$
est surjective, donc  $^t\pi$ est injective. 
Ainsi, $d= \dim(\mc F^\perp)=\dim \Xi(T)-\dim\mc F$ ; 
et l'{\'e}galit{\'e} du lemme est d{\'e}montr{\'e}e.

Comme $\dim(<\mc C>\cap(\mc F\times\Xi(\h T)_{\Q}))\geq 
\dim \mc C_{\mc F}$, on en d{\'e}duit que $\delta_{\mc F}\geq 0$.
\end{preuve}\\

On caract{\'e}rise alors le cas $\delta_{\mc F}=0$ dans la 

\begin{prop}\label{prop:critdelta=0}
On a {\'e}quivalence entre~:
\begin{enumerate}
\item $<\mc C_{\mc F}>=<\mc C>\cap\left(\mc
F\times\Xi(\hat T)_\Q\right)$ ;
\item $\delta_{\mc F}=0$ 
\end{enumerate}
On dit alors que la face $\mc F$ est pleine.
\end{prop} 

\begin{preuve}
Il est clair que $\mc C_{\mc F}$ et $<\mc C_{\mc F}>$ sont inclus dans 
$<\mc C>\cap\left(\mc F\times\Xi(\hat T)_\Q\right)$. 
Alors, l'assertion $(i)$ est {\'e}quivalente {\`a}
$\dim <\mc C_{\mc F}>=
\dim \left (<\mc C>\cap (\mc F\times\Xi(\hat T)_\Q)\right)$. 
D'apr{\`e}s le lemme~\ref{dimClin} ceci {\'e}quivaut {\`a} 
$\dim \mc C_{\mc F} =\dim\mc C-\dim\Xi(T)+\dim\mc F$, 
soit $\delta_{\mc F}=0$.
\end{preuve}

\subsection{Relations avec le polytope moment}\label{polmoment}

Soit $\hat\nu\in \hat{\mc D}$, posons~:
$$P_{\hat\nu}=\{\frac{\mu}{n}\in\Xi(T)_{\Q},\,(\mu,n\hat\nu)\in\mc
C\}\,.$$
Consid{\'e}rons $\hat{\mc B}:=\hat G/\hat B_-$ la vari{\'e}t{\'e} des drapeaux
de $\hat G$. Alors il existe un unique fibr{\'e} en droite $\hat
G$-lin{\'e}aris{\'e} $\mc L_{\hat\nu}$ sur $\hat{\mc B}$ tel que $\hat B_-$
agisse par le caract{\`e}re $-\hat\nu$ sur la fibre au dessus de $\hat
B_-/\hat B_-$. Pour tout entier $n$, le $\h G$-module des sections globales de 
$\mc L_{n\hat\nu}$ est isomorphe {\`a} $V_{n\hat\nu}$.
Mais alors $ P_{\hat\nu}$ est le polytope moment associ{\'e} {\`a} $\mc
L_{\hat\nu}$ pour l'action de $G$ sur $\hat{\mc B}$. Ce polytope moment est 
directement li{\'e} au polytope moment des g{\'e}om{\`e}tres symplecticiens, voir l'appendice 
de Mumford dans \cite{ne2}, et {\'e}galement \cite{Br:imagemom}. 
Les polytopes 
moments associ{\'e}s au d{\'e}composition de repr{\'e}sentations ont {\'e}t{\'e} {\'e}tudi{\'e}s notamment 
dans \cite{sjaamar}, \cite{man} et \cite{Br:genface}.  

La propri{\'e}t{\'e} d'{\^e}tre pleine s'interpr{\`e}te aussi en termes de ces
polytopes, si on suppose que $<\mc C>=\Xi(T)_\Q\otimes \Xi(\h T)_\Q$. 
Remarquons que cette hypoth\`ese est v\'erifi\'ee notamment lorsque $G$ est 
simple et non distingu\'e dans $\h G$, voir le corollaire \ref{prop:ker=0} 
ci-apr{\`e}s.
 
\begin{prop}\label{prop:polmom}
Supposons que $<\mc C>=\Xi(T)_\Q\otimes \Xi(\h T)_\Q$ ; si $\mc F$ est une face 
pleine $\mc D$, alors il existe un convexe $\Omega\subset\hat{\mc D}$ de 
dimension $\dim\h{\mc D}$
 tel que pour tout $\hat\nu\in \Omega$, $\dim P_{\h \nu}-
\dim(P_{\hat\nu}\cap\mc F)=\dim\Xi(T)-\dim\mc F$.
\end{prop}

\begin{preuve}

Soit $\pi$ la projection de $\Xi(T)_\Q\times \Xi(\h T)_\Q$ sur $\Xi(\h T)_\Q$.
Posons pour tout $\h\nu\in\h{\mc D}$, $H_{\h\nu}=\Xi(T)_\Q\times\Q.\h\nu$,
$\tilde P_{\hat\nu}=\mc C\cap H_{\h\nu}$ et 
$\tilde{\mc F}=\mc F\times\Xi(\h T)_\Q$. Il est clair que~: 
$$
\dim P_{\hat\nu}-\dim(P_{\hat\nu}\cap \mc F)
=\dim \tilde P_{\hat\nu}-\dim(\tilde P_{\hat\nu}\cap 
\tilde{\mc F})\ ;$$ 
d'autre part, comme $\mc F$ est pleine, on a {\'e}galit{\'e}~:
$$\dim\Xi(T)-\dim\mc F=\dim \mc C-\dim(\mc C\cap\tilde{\mc F}).$$
Il suffit donc de montrer que~: 
$$\dim\mc C- 
\dim(\mc C\cap\tilde{\mc F})=
\dim\tilde P_{\hat\nu}-
\dim(\tilde P_{\hat\nu}\cap\tilde{\mc F})
$$
pour $\h\nu$ dans un sous-ensemble convexe de $\h{\mc D}$ de dimension maximale.
Des raisonnements \'el\'ementaires de g\'eom\'etrie convexe montre que~:
\begin{equation}\label{ineg1}
\dim \mc C \geq \dim \tilde P_{\hat\nu}+\dim \pi(\mc C)-1 
\end{equation}
pour tout 
$\h\nu\in
\h{\mc D}$, avec {\'e}galit{\'e} si $\h\nu\in X_{\mc C}$, 
o\`u $X_{\mc C}$ 
est l'ensemble convexe des points $\h\nu\in\h{\mc D}$ tels que $H_{\h\nu}$ 
rencontre l'int{\'e}rieur de $\mc C$. De m\^eme,

\begin{equation}\label{ineg2}
\dim(\mc C\cap\tilde{\mc F}) \geq 
\dim(\tilde P_{\hat\nu}\cap 
\tilde{\mc F})+
\dim\pi(\mc C\cap \tilde{\mc F})-1
\end{equation}
pour tout 
$\h\nu\in \h{\mc D}$, avec {\'e}galit{\'e} si $\h\nu\in X_{\CF}$, o\`u 
$X_{\CF}$ est l'ensemble convexe des points $\h\nu\in\h{\mc D}$ tels que 
$H_{\h\nu}$ rencontre l'int{\'e}rieur de $\mc C\cap\tilde{\mc F}$.
On a l'inclusion $X_{\CF}\subset X_{\mc C}$. En effet si 
$\h\nu\in X_{\CF}$, alors 
$\h\nu\in X_{\mc C}$, {\`a} moins que $\tilde P_{\hat\nu}\subset 
\tilde P_{\hat\nu}\cap \tilde{\mc F}$. 
Cette inclusion implique :
$$
<\mc C>\cap H_{\h\nu}\subset <\mc C>\cap H_{\h\nu}\cap
\tilde{\mc F}
$$
ce qui contredit $<\mc C>=\Xi(T)_\Q\otimes \Xi(\h T)_\Q$.

Soit $\h\nu\in X_{\CF}$, alors on est  dans  le cas d'{\'e}galit{\'e} des 
in{\'e}quations \ref{ineg1} et \ref{ineg2}. On en d{\'e}duit que~:
$$
\dim\mc C- 
\dim(\mc C\cap\tilde{\mc F})=
\dim\tilde P_{\hat\nu}-
\dim(\tilde P_{\hat\nu}\cap\tilde{\mc F})
+\dim\pi(\mc C)-\dim\pi(\mc C\cap\tilde{\mc F}
).
$$

Puisque $\mc F$ est pleine, d'apr{\`e}s la proposition \ref{prop:critdelta=0}, 
$$
<\mc C\cap\tilde{\mc F}>=
<\mc C>\cap\tilde{\mc F}$$
et donc~:
$$
\dim\pi(\mc C\cap\tilde{\mc F})=\dim\pi(<\mc C>\cap \tilde{\mc F})
= \dim\pi(\tilde{\mc F})
=\dim\pi(\mc C)=\dim\h{\mc D}.
$$
Ce qui conclut la preuve.
\end{preuve}

\medskip

\begin{rem}
Si l'on omet l'hypoth{\`e}se $<\mc C>=\Xi(T)_\Q\otimes \Xi(\h T)_\Q$, le r{\'e}sultat de 
la 
proposition n'est plus vrai. On peut consid{\'e}rer l'exemple $G=H_1\times H_2$ et 
$\h G=H_1\times\h H_2$, o{\`u} $H_2$ est un sous-groupe de $\h H_2$ et $\mc F$ une 
face de 
la chambre dominante de $H_1$.
\end{rem}

\section{Noyaux d'action d'un tore}
\label{sec:noyau}

Consid{\'e}rons le sous-groupe de Levi $L$ de $G$ contenant $T$ et dont les 
racines sont celles de $G$ orthogonales {\`a} ${\mc F}$.
Soit $U$ le radical unipotent de $B$ et $P$ le sous-groupe parabolique de 
$G$ engendr{\'e} par $U$ et $L$.
Posons enfin $D:=[L,L]$ et $T_D:=D\cap T$. 
Soit ${\h U}^-$ le radical unipotent du sous-groupe de Borel ${\h B}^-$ 
oppos{\'e} {\`a} $\h B$ et contenant $\h T$.
On montre alors la 

\begin{prop}\label{poideC_F}
L'ensemble des poids de l'action de $T\times\hat T$ sur 
$k[\hat G]^{[P,P]\times\hat U^-}$ est {\'e}gal {\`a} $\mc{C}_\mc F$.
\end{prop}

\begin{preuve}
D'apr{\`e}s le th{\'e}or{\`e}me de Frobenius, le $\hat G\times\hat G$-module
rationnel $k[\hat G]$ se d{\'e}compose comme suit~:

$$k[\hat G]=\bigoplus_{\hat\nu\in\hat{\mc D}} V_{\hat\nu}\otimes
V_{\hat\nu}^*\,,$$
d'o{\`u}~: 
$$k[\hat G]^{[P,P]\times\hat U^-}=\bigoplus_{\hat\nu\in
\hat{\mc D}}V_{\hat\nu}^{[P,P]}\otimes V_{\hat\nu}^{*\hat U^-}\,.$$
or $\hat T$ agit sur $V_{\hat\nu}^{*\hat U^-}$ par $\hat\nu$. Ainsi,
l'ensemble des poids de $T\times\hat T$ dans $k[\hat
G]^{[P,P]\times\hat U^-}$ est {\'e}gal {\`a} l'ensemble des couples
$(\mu,\hat\nu)\in\Xi(T)\times\Xi(\hat T)$ tels que $\hat\nu\in\hat D$ et
$\mu$ est un poids de $T$ dans $V_{\hat\nu}^{[P,P]}$. Or il est bien
connu que $V_{\mu}^{[P,P]}\neq 0$ si et seulement si $\mu$
appartient {\`a} $\mc F$ ce qui conclut la preuve.
\end{preuve}

\bigskip

Posons $\mc C_{\mc F}^\vee:=\{(t,\hat t)\in T\times\hat T\ : \ 
\forall(\mu,\nu)\in\mc C_{\mc F}\ \mu(t)\hat\nu(\hat t)=1\}$ ; alors on
a la

\begin{prop}\label{C_F_orth}
Le groupe $\mc C_{\mc F}^\vee$ est le noyau de l'action de $T\times\hat T$
sur la vari{\'e}t{\'e} $[P,P]\DQG\hat G\DQD\hat U^-$.
\end{prop}

\begin{preuve}
D'apr{\`e}s la proposition~\ref{poideC_F}, $\mc C_{\mc F}^\vee$ est l'ensemble
des couples $(t,\hat t)\in T\times\hat T$ qui agissent trivialement sur 
$k\left[[P,P]\DQG\hat G\DQD\hat U^-\right]$, d'o{\`u} la proposition.
\end{preuve}

\section{Choix de $\hat T$ et $\h B$}\label{choix}

La proposition~\ref{C_F_orth} est vrai quelque soit le choix de 
$\h T$ et $\h B$.
Dans cette section, nous allons choisir deux tels sous-groupes de 
$\h G$ de sorte qu'il soit facile gr{\^a}ce {\`a} la d{\'e}composition de Bruhat de
d{\'e}crire un ouvert de $[P,P]\DQG\hat G\DQD\hat U^-$ stable
par $T\times\h T$.\\

Un de nos objectifs est le corollaire~\ref{cor:monotonie} ci-apr{\`e}s
qui compare 
$\delta_{\mc F_1}$ et $\delta_{\mc F_2}$ si $\mc F_1$ et $\mc F_2$ sont 
deux faces de $\mc D$ telles que $\mc F_1\subset\mc F_2$.
C'est pourquoi nous nous donnons ici deux telles faces de $\mc D$.\\

Rappelons que nous avons fix{\'e} les sous-groupes $T\subset B$ et
les faces $\mc F_1\subset\mc F_2$. 
Pour $i=1$, $2$, nous d{\'e}finissons alors comme dans la
section~\ref{sec:noyau}
les sous-groupes $T_{D_i}$, $D_i$, $L_i$ et $P_i$.
L'inclusion de $\mc F_1$ dans $\mc F_2$ implique que 
$P_2\subset P_1$ et $L_2\subset L_1$. 
En revanche, $\h T$ et $\h B$ ne sont pas suppos{\'e}s donn{\'e}s ici.

Soit $i=1$ ou $2$. 
Notons $S_i$ le centre connexe de $L_i$.
Consid{\'e}rons alors le centralisateur $\h L_i$ de $S_i$ dans $\h G$.
Notons $\h D_i$ le sous-groupe d{\'e}riv{\'e} de $\h L_i$ et $\h S_i$ le centre connexe
de $\h L_i$. On a alors~:
\begin{eqnarray}
\label{eq:LDS}
L_i=\h L_i\cap G\\
\h L_2\subset \h L_1
\end{eqnarray}

On veut maintenant construire des sous-groupes paraboliques 
$\h P_i$ de $\h G$ ayant des propri{\'e}t{\'e}s analogues {\`a} (\ref{eq:LDS}).

Pour $i=1$, $2$, on note St$_i(G)$ (resp. St$_i(\h G)$) l'ensemble des
poids non triviaux de $S_i$ dans Lie$(G)$ (resp. Lie$(\h G)$).
Pour tout $\alpha\in \Xi(S_i)$, on note $\mc H_i^\alpha$ le sous-espace
vectoriel de $\Xi_*(S_i)_\Q$ engendr{\'e} par les sous-groupes {\`a} un param{\`e}tre 
$\lambda$ de $S_i$ tels que $\alpha\circ\lambda$ est trivial.
Posons~:
$$
C^i=\{
\lambda\in\Xi_*(S_i)\,:\,
P_i=\left\{
g\in G\,:\,\lim_{t\to 0} \lambda(t)g\lambda(t)^{-1} {\rm\ existe}
\}
\right\}.
$$
Alors, $C^i$ est l'intersection de $\Xi_*(S_i)$ et d'une composante connexe 
$C^i_\Q$ du compl{\'e}mentaire  de $\cup_{\alpha\in{\rm St}_i(G)}\mc
H_i^\alpha$ dans  $\Xi_*(S_i)_\Q$.
De plus, comme $P_2\subset P_1$, via l'inclusion naturelle de 
$\Xi_*(S_1)$ dans $\Xi_*(S_2)$, 
$C^1_\Q$ est inclus dans l'adh{\'e}rence de $C^2_\Q$.

Consid{\'e}rons $\h C^1_\Q$ et $\h C^2_\Q$ deux composantes connexes de
$\Xi_*(S_1)_\Q\setminus\cup_{\alpha\in{\rm St}_1(\h G)}\mc H_1^\alpha$
et
$\Xi_*(S_2)_\Q\setminus\cup_{\alpha\in{\rm St}_2(\h G)}\mc H_2^\alpha$
telles que
$\h C^1_\Q\subset C^1_\Q$, $\h C^2_\Q\subset C^2_\Q$ et 
$\h C^1_\Q$ soit inclus dans l'adh{\'e}rence de $\h C^2_\Q$.
Pour $i=1,2$, on fixe $\lambda_i$ dans $\Xi_*(S_i)\cap \h C^i_\Q$.
Posons alors,
$$
\h P_i=\{g\in \h G\,:\,\lim_{t\to 0}\lambda_i(t)g\lambda_i(t^{-1})\ 
\mathrm{existe}\}.
$$
On a alors, 
\begin{eqnarray}
\label{eq:P1}
\h P_2\subset\h P_1\\
\label{eq:P2}
P_i=\h P_i\cap G.
\end{eqnarray}
De plus, comme 
$\h P_i^u=\{g\in \h G\ : \ \lim_{t\to 0}\lambda_i(t)g\lambda_i(t^{-1})=1\}$,
on a~:
\begin{eqnarray}
\label{eq:Pu}
P^u_i=\h P^u_i\cap G.
\end{eqnarray}

Soit enfin $\h B_{D_1}$ un sous-groupe de Borel de $\h D_1$ contenant $B_{D_1}$.
Posons $\h B_{D_2}=\h B_{D_1}\cap\h D_2$. 
Alors, $\h B_{D_2}$ est un sous-groupe de Borel de $\h D_2$.
De plus, on a~:

\begin{eqnarray}
\label{eq:B}
\h B_{D_2}=\h B_{D_1}\cap\h D_2,\\
B_{D_i}=\h B_{D_i}\cap G.
\end{eqnarray}

Soit $\h T_{D_1}$ un tore maximal de $\h B_{D_1}$ contenant $T_{D_1}$.
Posons alors $\h T_{D_2}=\h T_{D_1}\cap\h D_2$ et
$\h T=\h S_1\h T_{D_1}=\h S_2\h T_{D_2}$.
Alors, on a 
\begin{eqnarray}
\label{eq:T}
T_{D_i}=(\h T_{D_i}\cap G)^\circ.
\end{eqnarray}

\section{R{\'e}ductions}\label{reductions}
Nous allons dans cette section exprimer $\CFun^\vee$ 
et $\CFdeux^\vee$ comme noyaux d'une action  du tore $T\times\h T$  
sur deux vari{\'e}t{\'e}s \og comparables\fg\ ce qui permettra de 
comparer les dimensions de $\CFun$ et $\CFdeux$. 
Nous rappelons d'abord plusieurs lemmes bien connus sur les actions de groupes 
alg{\'e}briques. 

\subsection{Rappels}\label{rappels}

\begin{lem}\label{quot_ratio}
Soit $G$ un groupe alg{\'e}brique affine agissant sur une vari{\'e}t{\'e}
affine $X$ tels que $k[X]^G$ soit de type fini. 
Notons $\pi\,:\,X\rightarrow X\DQD G$
l'application quotient. 
Alors, se valent~:
\begin{enumerate}
\item $k(X)^G=\mathrm{Frac}(k[X]^G)$ ;
\item il existe un ouvert non vide $\ouv$ de $X\DQD G$ tel que pour
  $x\in\ouv$, $\pi^{-1}(x)$ contient une unique orbite ouverte de $G$. 
\end{enumerate}

On dit alors que le quotient $X\DQD G$ est rationnel.
\end{lem} 

\begin{preuve}
Voir  la section 2.4 p156 de \cite{PV}.
\end{preuve}

\begin{lem}
Soit $G$ un groupe alg{\'e}brique affine tel que $\Xi(G)=1$ et soit
$X$ une $G$-vari{\'e}t{\'e} affine. Si $k[X]$ est factorielle alors
$k[X]^G$ l'est aussi.
\end{lem}

\begin{preuve}
Voir \cite{PV} Theorem 3.17 p176.
\end{preuve}

\begin{lem}\label{frac.inv}
Soit $G$ un groupe alg{\'e}brique affine tel que $\Xi(G)=1$ et soit
$X$ une $G$-vari{\'e}t{\'e} factorielle affine.
Alors, le corps ${k(X)}^G$ des fractions rationnelles $G$-invariantes sur $X$
est {\'e}gale au corps des fractions de ${k[X]}^G$.
\end{lem}

\begin{preuve}
Voir \cite{PV} Theorem 3.3 p165.
\end{preuve}\\

Le lemme suivant est une cons{\'e}quence directe des lemmes pr{\'e}c{\'e}dents.

\begin{lem}\label{biration}
Soit $G$ un groupe alg{\'e}brique tel que $\Xi(G)=1$. 
Soit $X$ et $Y$ deux $G$-vari{\'e}t{\'e}s affines telles que~:
\begin{enumerate}
\item \label{biration1} les alg{\`e}bres $k[X]^G$ et $k[Y]^G$ sont de type 
fini ;
\item \label{biration2} le quotient $Y\DQD G$ est rationnel~;
\item \label{biration3} 
il existe une application $\varphi\,:\,X\rightarrow Y$, 
$G$-{\'e}quivariante et birationnelle ;
\end{enumerate}
Alors les vari{\'e}t{\'e}s $X\DQD G$ et $Y\DQD G$ sont birationnelles.
\end{lem}

Enfin, nous terminerons par une lemme concernant l'action d'un tore
sur une vari{\'e}t{\'e} affine dont la d{\'e}monstration est {\'e}vidente. 

\begin{lem}\label{kertore}
Soit $T$ un tore agissant sur deux vari{\'e}t{\'e}s $X$ et $Y$, alors on a ~:
\begin{enumerate}
\item il existe un ouvert non vide $\ouv$ de $X$ tel que pour tout
  $x\in\ouv$, l'isotropie de $T$ en $x$ est {\'e}gale au noyau de
  l'action de $T$ sur $X$ ;
\item si $f$ est un morphisme $T$-{\'e}quivariant de $X$ dans $Y$
  dont les fibres sont g{\'e}n{\'e}riquement finies, alors 
  $\Ker{T}{X}^\circ=\Ker{T}{Y}^\circ$. 
\end{enumerate}
\end{lem}

Pour $i=1$, $2$, on d{\'e}signe par $T_{D_i}\subset D_i\subset L_i\subset P_i$ 
les sous-groupes associ{\'e}s {\`a} $\mc F_i$ comme dans la 
section~\ref{sec:noyau}.
On a alors les inclusions suivantes~: $P_2\subset P_1$, $L_2\subset L_1$, 
$D_2\subset D_1$ et $P_1^u\subset P_2^u$. 
On choisit alors $\h L_1,\h L_2,\h P_1,\h P_2$\dots 
comme dans la section~\ref{choix}. 
Nous noterons $\gp_1^u:=\mathrm{Lie}(P_1^u)$ et $\h\gp_1^u:=\mathrm{Lie}(\h P_1^u)$. 
Remarquons que~:
\begin{itemize}
\item le groupe $D_1$ agit sur $\hat P_1^u$ par conjugaison en
laissant stable $P_1^u$ ; il agit donc sur $\h P_1^u/P_1^u$ ; 
\item la vari{\'e}t{\'e} $\h P_1^u/P_1^u$ est $D_1$-isomorphe {\`a} l'espace 
affine $\h\gp_1^u/\gp_1^u$ voir \cite{plm} ; 
\item $D_1$ agit {\'e}galement par multiplication {\`a} gauche sur $\hat L_1\DQD\hat U^-_{D_1}$ donc sur le produit $\h\gp_1^u/\gp_1^u\times\hat L_1\DQD\hat
U^-_{D_1}$. 
\end{itemize}

\begin{prop}\label{prop:modele1}
On a l'{\'e}galit{\'e}~:
$$\CFdeux^\vee=\Ker{T\times\hat T}{(\h\gp_1^u/\gp_1^u\times\hat L_1\DQD\hat U^-_{D_1})
\DQD [D_1\cap P_2,D_1\cap P_2]}.$$
En particulier (pour $\mc F_1=\mc F_2=\mc F$), 
$$\CF^\vee=\Ker{T\times\hat T}{\left( \h\gp^u/\gp^u\times\hat L\DQD\hat U^-_{D}\right ) \DQD D}\;;$$
Enfin, pour $\mc F=\mc F_1$ et $\mc F_2=\Xi(T)_\Q$, on obtient~: 
$$\mc C^\vee= \Ker{T\times\hat T}{\left(\h\gp^u/\gp^u\times\hat L\DQD\hat U^-_{D}\right )
\DQD U_D}.$$

\end{prop}

\begin{preuve} 
Rappelons que d'apr{\`e}s la proposition~\ref{C_F_orth} 
$$
\CFdeux^\vee=\Ker{T\times \hat T}{[P_2,P_2]\DQG\hat G\DQD\hat U^-}.
$$
La d{\'e}composition de Bruhat de $\h G$ par rapport au sous-groupe
parabolique $\h P_1$ implique que 
$\h P_1\h {P^u_1}^-=\h P^u_1\h {P_1}^-=\h P^u_1\h L_1\h {P^u_1}^-$ est un ouvert 
de $\h G$ stable par $P_2\times \h B^-$.
Mais alors, les lemmes~\ref{biration}
et~\ref{kertore} montrent que
$$
\CFdeux^\vee=
\Ker{T\times \hat T}{[P_2,P_2]\DQG \h P^u_1\h L_1\h {P^u_1}^- \DQD\hat U^-}.
$$
Or, le produit dans $\h G$ induit les trois isomorphismes suivants~:
$\hat U^-\simeq \h P_1^{u-}\psd \hat U^-_{D_1}$, 
$[P_2,P_2]\simeq P^u_1\psd [D_1\cap P_2,D_1\cap P_2]$ et
$\h P^u_1\h L_1\h {P^u_1}^-\simeq P^u_1\times \h L_1 \times \h {P^u_1}^-$.
Mais alors, 
$$\CFdeux^\vee=
\Ker{T\times \hat T}{\Big ((P_1^u\backslash \h P^u_1)\times (\h L_1\DQD\hat U^-_{D_1}\Big ))
\DQD [D_1\cap P_2,D_1\cap P_2]}.
$$
\end{preuve}
\medskip

Soit $\textgoth u$ et $\h{\textgoth u}$ les alg{\`e}bres de Lie de $U$ et $\h U$. Le corollaire 
suivant qui se d{\'e}duit d'un cas particulier de la proposition \ref{prop:modele1} permet 
de calculer la dimension de $\mc C$.

\begin{cor}\label{prop:ker=0}
On a les {\'e}galit{\'e}s~: 
$$\dim\mc C^\vee=\dim\Ker{T}{\h{\textgoth u}/\textgoth u}=\dim\Ker{T}{\h G/G},$$
o\`u $T$ agit par conjuguaison sur $\h{\textgoth u}/\textgoth u$ et par multiplication 
{\`a} gauche sur $\h G/G$. 
En particulier si $G$ est simple et non distingu{\'e} dans $\h G$, alors 
$\dim\mc C^\vee=0$.
\end{cor}

\begin{preuve}
La premi{\`e}re {\'e}galit{\'e} se d{\'e}duit directement de la proposition \ref{prop:modele1}
appliqu{\'e}e {\`a} $\mc F_1=\mc F_2=\Xi(T)_\Q$. 

 Pour la deuxi{\`e}me {\'e}galit{\'e} on remarque tout d'abord que~: 
$$\Ker{T}{\h{\textgoth u}/\textgoth u}=\Ker{T}{\h{\textgoth g}/\textgoth g},$$ 
o{\`u} 
$\h{\textgoth g}$, $\textgoth g$ sont les alg{\`e}bres de Lie respectives de $\h G$ et $G$ et 
o{\`u} l'action de $T$ sur $\h{\textgoth g}/\textgoth g$ est induite par l'action adjointe.
D'apr{\`e}s Luna \cite{lun:fermee}, il existe un morphisme {\'e}tale $T$-{\'e}quivariant de 
$\h G/G$ dans 
$\h{\textgoth g}/\textgoth g$, et donc d'apr{\`e}s le lemme \ref{kertore}, on a {\'e}galit{\'e}~:
$$\dim\Ker{T}{\h{\textgoth g}/\textgoth g}=\dim\Ker{T}{\h G/G},$$ 
o{\`u} $T$ agit sur la vari{\'e}t{\'e} ${\h G/G}$
par conjuguaison. Mais cette action est {\'e}gale {\`a} l'action de $T$ sur ${\h G/G}$
par multiplication {\`a} gauche.

Pour montrer la derni{\`e}re assertion, il suffit de remarquer que le noyau de 
$T$ sur $\h G/G$ est inclus dans  
$\bigcap_{\h g\in\h G}\h gG\h g^{-1}$ qui est un sous-groupe distingu{\'e} de $G$ et 
de $\h G$.
\end{preuve}

\subsection{}
Consid{\'e}rons les trois applications finies~:

$$
\begin{array}{cccc}
\tau\,:&S_1\times (S_2\cap T_{D_1})^\circ\times T_{D_2}&\longto&T\ ,\\
&(s,t,u)&\longmapsto&stu
\end{array}
\hspace{22pt}
\begin{array}{cccc}
\hat{\tau}\,:&\hat{S}_1\times \hat{T}_{D_1}&\longto&\hat{T}\\
&(s,t)&\longmapsto&st
\end{array}
$$

et 

$$
\begin{array}{cccc}
\varphi\,:&\hat P^u_1\times \hat{D}_1\times\hat{S}_1&\longto&
\hat{P}^u_1\times\hat{L}_1\\
&(\hat p,\hat d,\hat s)&\longmapsto&(\hat p,\hat d\hat s)
\end{array}
$$

On munit $\hat{P}^u_1\times\hat{D}_1\times \hat{S}_1$ 
d'une action de $\hat U^-_{D_1}\times (\hat P^u_1D_1)$ par~:

$$(\hat u,pd).(\hat p,\hat d,\hat s)=
(pd\hat pd^{-1},d\hat d\hat{u}^{-1},\hat s).
$$

On munit enfin $\hat P^u_1\times\hat D_1\times\hat S_1$ d'une action
du tore $S_1\times (S_2\cap T_{D_1})^\circ\times T_{D_2}\times
\hat{S}_1\times \hat{T}_{D_1}$ (not{\'e} $\mathbb T$) par~:
$$
(s_1,s_2,t,\hat s_1,\hat t_1).(\hat p,\hat d,\hat s)=
(s_1s_2t \hat p t^{-1}s_2^{-1}s_1^{-1}, s_2t_1 \hat d\hat t_1^{-1},s_1\hat s \hat s_1^{-1}).
$$

On peut maintenant simplifier l'expression de la composante neutre de
$\CFdeux^\vee$~:

\begin{prop}
\label{prop:modele2}
Posons $H_2=[D_1\cap P_2,D_1\cap P_2]$.
La composante neutre de $\CFdeux^\vee$ est isomorphe {\`a}
$$
T_{D_2}\times
\Ker{S_1\times (S_2\cap T_{D_1})^\circ\times\hat T_{D_1}}{
\left(\h\gp_1^u/\gp_1^u\times\hat D_1\DQD\hat U^-_{D_1}
\right)
\DQD H_2}^\circ.
$$
\end{prop}

\begin{preuve}
Notons $\CFdeux^{\vee,\circ}$ la composante neutre de $\CFdeux^\vee$.
On munit les vari{\'e}t{\'e}s $\h\gp_1^u/\gp_1^u\times\hat L_1\DQD\hat U^-_{D_1}$ 
et $\h\gp_1^u/\gp_1^u\times\hat D_1\DQD\hat U^-_{D_1}\times S_1$
d'une action du tore $\mathbb{T}$ gr{\^a}ce {\`a} $\tau$ et $\hat\tau$ et par passage au quotient. 
L'application $\varphi$ induit alors une application {\'e}quivariante
dont les fibres sont g{\'e}n{\'e}riquement finies entre ces deux $\mathbb{T}$-vari{\'e}t{\'e}s.
Mais alors, la proposition~\ref{prop:modele1} et le lemme~\ref{kertore} montrent que
les composantes neutres de $\CFdeux^\vee$ et de 
$
\Ker{\mathbb{T}}{(\h\gp_1^u/\gp_1^u\times\hat D_1\DQD\hat U^-_{D_1}\times S_1)\DQD H_2}
$ sont isomorphes. 

Posons $M=\h\gp_1^u/\gp_1^u\times\hat D_1\DQD\hat U^-_{D_1}$.
Comme $H_2$ est inclus dans $D_1$, 
$(\h\gp_1^u/\gp_1^u\times\hat D_1\DQD\hat U^-_{D_1}\times S_1)\DQD H_2$ est isomorphe
{\`a} $M\DQD H_2\times \hat S_1$.
Comme $T_{D_2}$ est inclus dans $H_2$, on en d{\'e}duit que 
$$
\CFdeux^{\vee,\circ}\simeq
\left(
T_{D_2}\times 
\Ker{S_1\times (S_2\cap T_{D_1})^\circ\times \hat{S}_1\times \hat{T}_{D_1}}
{M\DQD H_2\times \hat S_1}
\right)^\circ,
$$
puis que 
$$
\CFdeux^{\vee,\circ}\simeq
T_{D_2}\times \left(
\Ker{S_1\times (S_2\cap T_{D_1})^\circ\times \hat{T}_{D_1}}
{M\DQD H_2}
\right )^\circ
\,.$$
\end{preuve}\\

Pour arriver jusqu'{\`a} notre th{\'e}or{\`e}me principal, il nous faut
maintenant d{\'e}crire plus pr{\'e}cis{\'e}ment le noyau apparaissant dans
la proposition~\ref{prop:modele2}.
Pour s'all{\'e}ger de notations inutiles, 
dans la section~\ref{annexe}, nous pr{\'e}senterons dans un cadre 
plus g{\'e}n{\'e}ral  les r{\'e}sultats qui nous permettent de conclure.

\begin{theo}
\label{th:ppal}
Il existe un ouvert non vide  $\ouv$ de 
$\hat{\gp}^u_1/\gp^u_1\times \hat{D}_1/\hat{B}_{D_1}$ tel que pour tout $x$ 
dans $\ouv$, ${\mc C}^{\vee,\circ}_{\mc F_2}$ est isomorphe au produit de 
$T_{D_2}$ par la composante neutre de l'isotropie r{\'e}ductive 
en $x$ du groupe $L_1\cap P_2$.
\end{theo}
\begin{preuve} 
  On applique la proposition~\ref{prop:rho1} {\`a} la vari{\'e}t{\'e} 
${\bf M}=\hat{\gp}^u_1/\gp^u_1\times \hat{D}_1\DQD U_{\h D_1}^{-}$, au tore
${\bf T}=S_1\times\ T_{D_1}$, au groupe semi-simple
${\bf D}=D_1$ et au sous-groupe
parabolique ${\bf P}=P_2\cap D_1$. Remarquons qu'alors $L_2\cap
D_1$ est un Levi de ${\bf P}$ et que son centre connexe est {\'e}gal
{\`a} ${\bf S}=S_2$.
On obtient l'existence d'un ouvert non vide
$\ouv$ de $\hat{\gp}^u_1/\gp^u_1\times \hat{D}_1\DQD U_{\h D_1}^{-}$ tel
que pour tout $x\in\ouv$~:
$$
\Ker{S_1\times (S_2\cap T_{D_1})^\circ\times \hat{T}_{D_1}}
{M\DQD H_2}=\rho(((P_2\cap D_1)\times S_1\times\ T_{D_1})_x).
$$
Remarquons que puisque $\hat{D}_1\DQD U_{\h D_1}^{-}$ contient 
$\hat{D}_1/U_{\h D_1}^{-}$ comme ouvert, quitte {\`a} r{\'e}tr{\'e}cir $\ouv$,
on peut supposer que l'{\'e}galit{\'e} ci-dessus est vrai pour tout
$x\in\ouv\subset\hat{\gp}^u_1/\gp^u_1\times \hat{D}_1/U_{\h D_1}^{-}$.
 On applique ensuite la proposition~\ref{prop:rho2}, avec
${\bf D}\subset\h{\bf{D}}=\h D_1$ et
${\bf M}=\hat{\gp}^u_1/\gp^u_1\times \hat{D}_1/U_{\h D_1}$. On
obtient alors que pour $x=(a,b)\in \hat{\gp}^u_1/\gp^u_1\times
\hat{D}_1/U_{\h D_1}$, $\rho(((P_2\cap D_1)\times
S_1\times\ T_{D_1})_{(a,b)})$ est isomorphe au quotient de
$(P_2\cap D_1\times S_1)_{(a,q(b))}$ par son radical unipotent, o{\`u} 
$q$ est l'application quotient de $\hat{D}_1/U_{\h D_1}^{-}$ dans 
$\hat{D}_1/B_{\h D_1}^{-}$. Le th{\'e}or{\`e}me suit, en remarquant que 
$(P_2\cap D_1)\times S_1\simeq P_2\cap L_1$.
\end{preuve}

\begin{cor}
\label{cor:monotonie}
Avec les notations ci-dessus, 
$$
\delta_{\mc F_1}\geq \delta_{\mc F_2}
$$ 
\end{cor}

\begin{preuve}
Comme $\dim \mc F_1-\dim \mc F_2=\dim T_{D_2}-\dim T_{D_1}$, on a
$\delta_{\mc F_1}-\delta_{\mc F_2}=
(\dim \mc C^\vee_{\mc F_1}-\dim T_{D_1})-
(\dim \mc C^\vee _{\mc F_2}-\dim T_{D_2)}$.
Mais alors, le corollaire d{\'e}coule imm{\'e}diatement du 
th{\'e}or{\`e}me~\ref{th:ppal} appliqu{\'e} une fois {\`a} la paire de face 
$(\mc F_1,\mc F_2)$ et une fois {\`a} la paire $(\mc F_1,\mc F_1)$.
\end{preuve}

\begin{cor}\label{cor:ppal}
Il existe un ouvert  non vide $\ouv$ de 
$\hat{\gp}^u/\gp^u\times \hat{D}/\hat{B}_{D}$ tel que pour tout $x$
dans $\ouv$, $\delta_{\mc F}$ est {\'e}gal {\`a} la diff{\'e}rence des
dimensions des isotropies r{\'e}ductives des groupes $L$ et $B_{L}$ en $x$. 
\end{cor}

\begin{preuve}
Il suffit d'appliquer le th{\'e}or{\`e}me~\ref{th:ppal} une premi{\`e}re fois {\`a}
la paire de faces $(\mc F,\mc F)$ et ensuite {\`a} la paire 
$(\mc F,\Xi(T)_\Q)$.
\end{preuve}\\

Sur les exemples, le th{\'e}or{\`e}me~\ref{th:ppal} sera souvent 
utiliser via le corollaire suivant~:

\begin{cor}
\label{cor:isotfini}
S'il existe un point de $\hat{\gp}^u/\gp^u\times \hat{D}/\hat{B}_{D}$ 
dont l'isotropie  
dans $D$ est finie, alors $\mc F$ est pleine.
En particulier, s'il existe un point de $\h G/\h B$ dont l'isotropie 
dans le groupe d{\'e}riv{\'e} de $G$ est finie, alors toutes les faces 
de $\mc D$ sont pleines.
\end{cor}

\begin{preuve}
Soit $\ouv$ un ouvert de $\hat{\gp}^u/\gp^u\times \hat{D}/\hat{B}_{D}$ 
v{\'e}rifiant le corollaire~\ref{cor:ppal}. 
Comme l'ensemble des points dont l'isotropie
pour $D$ est fini est un ouvert, il existe $x\in\ouv$ tel que $D_x$ soit fini.
Mais alors, $L_x^\circ$ est inclus dans le centre de $L$ et donc 
dans $B_L$ et $\mc F$ est pleine.
S'il existe $x\in\h G/\h B$ dont l'isotropie pour $D$ est finie, 
nous venons de montrer que la face $\{0\}$ est pleine. 
Mais alors, le corollaire~\ref{cor:monotonie} montre que toutes les 
faces sont pleines.
\end{preuve}

\section{Exemples}
\label{sec:exple}

\begin{exple} Le produit tensoriel.
\begin{prop}
\label{exple:GGG}
Soit $G$ un groupe simple.
Pour l'inclusion diagonale de $G$ dans $G\times G$ toutes les faces non 
r{\'e}duites {\`a} un point 
sont pleines.
De plus, $\delta_{\{0\}}$ est {\'e}gal au rang de $G$.
\end{prop}

\begin{dem}
Fixons un sous-groupe de Borel $B$ de $G$ et un tore maximal $T$ de $B$.
Soit $\cF$ une face de l'ensemble des poids dominants de $(G,B,T)$.
Notons $P$ le sous-groupe parabolique de $G$ associ{\'e} {\`a} la face
$\cF$, $L$ son sous-groupe de Levi contenant $T$ et enfin $D$ le sous-groupe
d{\'e}riv{\'e} de $L$. 
Posons~:
$$
\begin{array}
{l@{\ \ \ }l@{\ }l@{\> \ \ }l@{\ }l@{\> \ \ }l@{\ }l@{\> \ \ }l}
B_D=B\cap D&
\hP=P\times P&\hD=D\times D&\hB=B\times B\ &\h B_D=B_D\times B_D.
\end{array}
$$
Alors, $\h B_D\subset \hD\subset \hP\subset \hG$ satisfont aux
propri{\'e}t{\'e}s de la section~\ref{choix}. 
De plus, $M=\gp^u\times D/B_D\times D/B_D$.

En particulier, l'isotropie g{\'e}n{\'e}rique $I$ de $D$ agissant sur 
$M$ est {\'e}gale {\`a} celle de $T_D$ agissant sur $\gp_u$. 
Soit $P^-$ le sous-groupe parabolique contenant $T$ et oppos{\'e} {\`a} $P$. 
Comme $\gp_u$ est isomorphe {\`a} un ouvert $T$-stable de $G/P^-$,
$I$ est {\'e}gale au noyau de l'action de $T_D$ dans $G/P^-$, c'est-{\`a}-dire
{\`a} $\bigcap_{g\in G}gP^-g^{-1}\cap T_D$.
Comme $G$ est simple cette derni{\`e}re est finie {\`a} moins que $P^-=G$, c'est-{\`a}-dire
{\`a} moins que $\cF=\{0\}$.

Supposons maintenant que $\cF=\{0\}$. 
Alors, $M=G/B\times G/B$. Donc, l'isotropie g{\'e}n{\'e}rique de $D=G$ sur $M$ est un tore maximal 
de $G$ alors que celle de $B$ est finie. 
Le corollaire~\ref{cor:ppal} ach{\`e}ve alors la d{\'e}monstration.
\end{dem}
\end{exple}\\

\begin{exple}\label{expleS2}
Soit $V$ un $k$-espace vectoriel de dimension au moins deux.
Consid{\'e}rons le morphisme $i\,:\,{\rm Sl}(V)\rightarrow {\rm Sl}(S^2V)$. 
On a~:
\begin{prop}
Pour l'inclusion de $G=i({\rm Sl}(V))$ dans $\h G={\rm Sl}(S^2V)$, toutes les faces 
sont pleines. 
\end{prop}

D'apr{\`e}s le corollaire~\ref{cor:isotfini}, il suffit de montrer le 

\begin{lem}
\label{lem:S2}
Il existe un drapeau complet de $S^2V$ dont l'isotropie pour ${\rm Sl}(V)$ est finie.
\end{lem}

\begin{preuve} 
Nous allons construire un drapeau complet $\Omega$ de $S^2V^*$ dont 
le groupe d'isotropie $H$ est fini. Son orthogonal v{\'e}rifiera
la propri{\'e}t{\'e} du lemme. 
Soit $\mathcal B=(e_1,\cdots ,e_n)$ une base de $V$ et $(e_1^*,\cdots,e^*_n)$ 
sa base duale, le premier sous-espace du drapeau $\Omega$ est 
la droite engendr{\'e}e par $\omega_0=\Sigma {e_i^{*}}^2\in S^2V^*$. 
La composante neutre du groupe d'isotropie de cette droite est ${\rm SO}(V)$. 
Il suffit donc de montrer que l'intersection $H\cap {\rm SO}(V)$ est finie. 
Si $\omega$ est un {\'e}l{\'e}ment quelconque de $S^2V^*$, on d{\'e}finit~:
$$
\begin{array}{cccc}
\tilde\omega\,:\,&V&\rightarrow&V^* \\
&x&\mapsto&\omega(x,.)\ .\\
\end{array}
$$
Comme $\omega_0$ est non d{\'e}g{\'e}n{\'e}r{\'e}e, $\tilde\omega_0$ est un 
isomorphisme ${\rm SO}(V)$-{\'e}quivariant.  
Ce qui nous permet de d{\'e}finir une application~:

$$
\begin{array}{cccc}
\varphi\,:\,&S^2 V^*&\rightarrow&{\rm End}(V) \\
&\omega&\mapsto&\tilde\omega_0^{-1}\circ\tilde\omega\ .\\
\end{array}
$$
L'application $\varphi$ est un isomorphisme lin{\'e}aire
${\rm SO}(V)$-{\'e}quivariant (l'action
sur l'espace de droite est donn{\'e}e par la conjugaison). Soit 
$\omega_1=\Sigma i{e^*_i}^2$. On voit alors que $\varphi(\omega_0)$ est
l'identit{\'e} de $V$ et
que $e_i$ est vecteur propre de 
$\varphi(\omega_1)$ pour la valeur propre $i$.
Le stabilisateur du drapeau $<\omega_0>\subset <\omega_0,\omega_1>$ dans ${\rm SO}(V)$
est alors compos{\'e} des {\'e}l{\'e}ments $g\in {\rm SO}(V)$ tel qu'il existe
$a,b\in k$ v{\'e}rifiant 
$g.\varphi(\omega_1)=\varphi(a\omega_0+b\omega_1)$. 
La matrice de $\varphi(a\omega_0+b\omega_1)$ dans $\mc B$ est une 
matrice diagonale de valeurs propres $\{a+bi\,:\,i=1,\cdots,n\}$. 
Cet ensemble est {\'e}gal {\`a} 
l'ensemble des valeurs propres de $g.\varphi(\omega_1)$ (donc de $\varphi(\omega_1)$ soit 
$\{1,\cdots,n\}$). On en d{\'e}duit que ou bien $a=0$ et $b=1$, ou bien
$a=1+n$ et $b=-1$ ; puis que 
la composante neutre du groupe d'isotropie de 
$<\omega_0>\subset <\omega_0,\omega_1>$ est inclus dans l'ensemble des matrices 
diagonales.
Comme par ailleurs elle est inclus dans ${\rm SO}(V)$, elle est triviale.
Le lemme suit alors facilement.
\end{preuve}\\
\end{exple}

\begin{exple}
Nous allons ici regarder l'exemple donn{\'e} par le morphisme 
$i\,:\,{\rm Sl}(V)\rightarrow {\rm Sl}(\extdeux V)$. 
On pose alors $G=i(\mathrm{Sl}(V)$ et 
$\h G=\mathrm{Sl}(\extdeux V)$.
Si $\dim V=2$, alors  $\extdeux V\simeq k$ et $\h G$
est trivial ; si $\dim V=3$, alors 
$\extdeux V\simeq V^*$ et $G=\h G$. On suppose donc que
$\dim V\geq 4$.

\begin{prop}
Si $\dim V\geq 4$, alors pour l'inclusion de $G=i({\rm Sl}(V))$ 
dans $\h G={\rm Sl}(\extdeux V)$, toutes les faces de $\mc D$ sont pleines.
\end{prop}

Comme dans l'exemple pr{\'e}c{\'e}dent, la proposition est une cons{\'e}quence du

\begin{lem}\label{isotropielambda2}
Si $\dim V\geq 4$, alors il existe un drapeau de $\extdeux(V^*)$
dont l'isotropie dans ${\rm Sl}(V)$ est finie.
\end{lem}

\begin{preuve}
Dans $\extdeux (V^*)$, il existe un {\'e}l{\'e}ment non d{\'e}g{\'e}n{\'e}r{\'e} 
si et seulement si $\dim V$ est paire. 
Commen{\c c}ons par le cas o{\`u} $\dim V$ est paire. 
On pose alors $\dim V=2p$ et on fixe une base  
$\mathcal B=(e_1,\cdots,e_p,\varepsilon_1,\cdots,\varepsilon_p)$ de
$V$ et 
$\mathcal B^*=
(e_1^*,\cdots,e_p^*,\varepsilon_1^*,\cdots,\varepsilon_p^*)$ sa base duale. 
On pose alors $\omega_0=\Sigma_{i=1}^pe_i^*\wedge\varepsilon_i^*$. 
Soit $H_0$ la composante neutre de l'isotropie de la droite engendr{\'e}e par
$\omega_0$. 
Comme $H_0\subset {\rm Sl}(V)$, on voit que $H_0={\rm Sp}(\omega_0)$. 
Comme dans la d{\'e}monstration du lemme~\ref{lem:S2}, on d{\'e}finit pour tout 
$\omega\in\extdeux V^*$ une application 
$\tilde\omega\,:\,V\longto V^*$, $x\longmapsto \omega(x,.)$.
On d{\'e}finit alors l'application lin{\'e}aire, ${\rm Sp}(\omega_0)$-{\'e}quivariante
$\varphi\,:\,\extdeux V^*\longto \End(V)$, 
$\omega\longmapsto\tilde\omega_0^{-1}\circ\tilde\omega$.

Soit $\omega_1=\Sigma_{i=1}^pie_i^*\wedge\varepsilon_i^*$ ; 
$\varphi(\omega_1)$ admet $i=1,\cdots,p$ comme valeurs propres,
et $(e_i^*, \varepsilon_i^*)$ engendre le sous-espace
propre associ{\'e} {\`a} $i$. 
Par un raisonnement analogue {\`a} l'exemple pr{\'e}c{\'e}dent, on
montre alors que la composante neutre du sous-groupe
d'isotropie du drapeau~: $<\omega_0>\subset <\omega_0,\omega_1>$ est
l'ensemble des matrices de ${\rm Sp}(V)$ pr{\'e}servant les sous-espaces
$W_i:=<e_i^*,\varepsilon_i^*>$. Cette composante neutre est donc isomorphe au produit~:
$\prod_{i=1}^p {\rm Sl}(W_i)$. Sous l'action de ce groupe, l'espace
$\extdeux V^*$ se d{\'e}compose en repr{\'e}sentations irr{\'e}ductibles~: 
$$
\extdeux V^*=\bigoplus_{i=1}^p\extdeux W_i\oplus\bigoplus_{1\leq
i<j\leq p}W_i\otimes W_j\ .
$$ 

La suite de la construction du drapeau se fait en consid{\'e}rant
dans l'espace dual $\extdeux (V^*)^*\simeq\extdeux V$ des
vecteurs $\omega_3^*,\omega_4^*,\cdots$ dont l'orthogonal dans $\extdeux
V^*$ contient $<\omega_0,\omega_1>$. {\`A} cause de la d{\'e}composition
ci-dessus, de tels vecteurs peuvent {\^e}tre pris dans 
$\bigoplus_{1\leq i<j\leq p}W^*_i\otimes W^*_j$. 
Gr{\^a}ce {\`a} cette d{\'e}composition, on peut
{\'e}galement supposer que $p=2$. Pour ${\rm Sl}(W)\times {\rm Sl}(W)$, l'espace
$W\otimes W$ est isomorphe {\`a} l'espace $\End(W)$. {\`A} indice finie
pr{\`e}s l'isotropie de la droite engendr{\'e}e par l'application
identit{\'e} est le groupe ${\rm Sl}(W)$ inclus diagonalement dans
${\rm Sl}(W)\times {\rm Sl}(W)$ et qui agit sur $\End(W)$ par conjugaison. 
Comme ${\rm Sl}(W)$-repr{\'e}sentation $\End(W)$ se d{\'e}compose comme la
somme de la repr{\'e}sentation triviale et de la repr{\'e}sentation
adjointe. Le probl{\`e}me se ram{\`e}ne donc {\`a} exhiber un drapeau de cette
repr{\'e}sentation adjointe d'isotropie finie. Soit $H$ et $F$ les
{\'e}l{\'e}ments suivants~: 

$$
H=\left(
\begin{array}{cc}
1&0\\
0&-1
\end{array}
\right)
F=\left(
\begin{array}{cc}
0&1\\
1&0
\end{array}
\right)
$$
(dans une base quelconque de W). Alors il est imm{\'e}diat que le
drapeau suivant~:$<H>\subset <H,F>$ 
a une isotropie finie dans $\mathrm{Sl}(W)$.

Supposons maintenant que la dimension de $V$ est impaire et {\'e}gale
{\`a} $2p+1$, avec $p\geq 2$. Soit $\mathcal
B=(e_1,\cdots,e_p,\varepsilon_1,\cdots,\varepsilon_p,\kappa)$ une base de
$V$. Nous noterons $W$ le sous-espace vectoriel de $V$ de base
$(e_1,\cdots,e_p,\varepsilon_1,\cdots,\varepsilon_p)$. 

Rappelons qu'on a une injection ${\rm Gl}(V)$-{\'e}quivariante 
de $\extdeux V^*$ dans $\Hom(V^*,V)$. Soient $\omega_0$, $\omega_1$ et
$\omega_2$ les {\'e}l{\'e}ments de $\Hom(V^*,V)$ dont les matrices dans
les bases $\mc B^*,\mc B$ sont~:

$$
\omega_0=\left(
\begin{array}{ccc}
0&I&0\\
-I&0&0\\
0&0&0
\end{array}
\right)
\omega_1=\left(
\begin{array}{ccc}
0&D&^tV\\
-D&0&0\\
-V&0&0
\end{array}
\right)
\omega_2=\left(
\begin{array}{ccc}
0&D'&0\\
-D'&0&^tV\\
0&-V&0
\end{array}
\right)
$$
o{\`u} $I$ d{\'e}signe la matrice identit{\'e} de taille $p\times p$, 
$D$ la matrice diagonale avec comme termes diagonaux~: 
$1,2,\cdots,p$, $V$ le vecteur
$(1,\cdots,1)\in k^p$ et $D'$ une matrice de $\mathrm M_p(k)$ telle que
les matrices $I,D,D'$ soient lin{\'e}airement ind{\'e}pendantes.
Nous allons maintenant d{\'e}finir un hyperplan de $\extdeux V^*$
contenant les vecteurs $\omega_0,\omega_1$ et $\omega_2$ ce qui
revient {\`a} d{\'e}finir une droite de $\extdeux V^{**}$, qui s'identifie {\`a}
$\extdeux V$ et donc {\`a} un sous-espace de $\Hom(V,V^*)$. 
Remarquons que la dualit{\'e} entre $\extdeux V$ 
et $\extdeux V^*$ est la restriction de l'application bilin{\'e}aire suivante~:
$$
\begin{array}{ccc}
\Hom(V,V^*)\times\Hom(V^*,V)&\longto &k\\
(A,B)&\longmapsto&\mathrm{trace}(B\circ A) 
\end{array}
$$
Soit $\theta\in\Hom(V^*,V)$ d{\'e}fini par~: 

$$
\left(
\begin{array}{ccc}
0&U&0\\
-U&0&0\\
0&0&0
\end{array}
\right)
$$
avec $U$ une matrice telle que
$\mathrm{tr}(U)=\mathrm{tr}(UD)=\mathrm{tr}(UD')=0$, (ce choix est
possible puisque $p\geq 2$ et donc $\dim\mathrm{M}_p(k)\geq 4$). On a
alors que $<\omega_0,\omega_1,\omega_2>\subset \ker \theta$.
Consid{\'e}rons donc le drapeau partiel~:
$$\Omega
=(<\omega_0>\subset<\omega_0,\omega_1>\subset<\omega_0,\omega_1,\omega_2>\subset
\ker \theta)\ .$$  

Comme $\ker\theta=k\kappa^*$ et $\ker\omega_0=k.\kappa$ le groupe
d'isotropie de ce drapeau est inclus
dans le groupe $\mathrm{Gl}(W)\times\mathrm{Gl}(k.\kappa)$. Ce groupe
respecte la d{\'e}composition $\extdeux V^*=\extdeux W^*\oplus
W^*\otimes k.\kappa^*$. Les {\'e}l{\'e}ments $\omega_0$ et $\omega_1$
s'{\'e}crivent dans $\extdeux V^*$~:
$$
\omega_0=\Sigma_{i=1}^pe_i^*\wedge\varepsilon^*_i\mathrm{\ \ et\ \  }
\omega_1=\Sigma_{i=1}^pie_i^*\wedge\varepsilon^*_i
+\kappa^*\wedge(e_1^*+\cdots+e^*_n)
$$
On en d{\'e}duit qu'un {\'e}l{\'e}ment du stabilisateur dans
$\mathrm{Gl}(W)\times\mathrm{Gl}(k\kappa)$ du drapeau~:
$<\omega_0>\subset<\omega_0,\omega_1>$ stabilise le drapeau de
$\extdeux W^*$~: 
$<\omega_0>\subset<\omega_0,\Sigma ie_i^*\wedge\varepsilon_i^*>$ et la
droite engendr{\'e}e par $\kappa^*\wedge(e_1^*+\cdots+e^*_n)$. En utilisant
cette propri{\'e}t{\'e} et le cas de la dimension paire, on en d{\'e}duit
que la composante neutre
du groupe d'isotropie du drapeau $\mc D$ est inclus dans le groupe
compos{\'e} des matrices de la forme~:
$$
\left(
\begin{array}{ccc}
\lambda^{p+1} I&0&0\\
Z&\lambda^{-p} I&0\\
0&0&\lambda^{-p}
\end{array}
\right)
$$
o{\`u} $\lambda\in k^*$ et $Z$ est une matrice de taille $p\times p$
diagonale. Ensuite on {\'e}crit
matriciellement que de tels {\'e}l{\'e}ments envoient $\omega_2$ dans 
$<\omega_0,\omega_1,\omega_2>$. En utilisant notamment le fait que les
matrices $I,D,D'$ sont lin{\'e}airement ind{\'e}pendantes, on d{\'e}duit que
la composante  neutre de l'isotropie du drapeau $\Omega$  est triviale. 
\end{preuve}\\
\end{exple}

\begin{exple}
Soient $E,F$ deux $k$-espaces vectoriels de dimensions finies~; le groupe
$\mathrm{Gl}(E)\times \mathrm{Gl}(F)$ agit sur $E\otimes F$. On a donc
un morphisme $i\,:\,\mathrm{Gl}(E)\times \mathrm{Gl}(F)\rightarrow
\mathrm{Gl}(E\otimes F)$. On pose $\h G=\mathrm{Gl}(E\otimes F)$ et
$G=i(\mathrm{Gl}(E)\times \mathrm{Gl}(F))$. On a alors la

\begin{prop}
Dans la situation ci-dessus, toutes les faces sont pleines.
\end{prop}

\begin{preuve}
La d{\'e}monstration qui est analogue aux deux exemples pr{\'e}c{\'e}dents est
ici omise.
\end{preuve}
\end{exple}

\section{R{\'e}sultats Annexes}\label{annexe}

\subsection{}
\label{sec:situation1}
Soit ${\bf T}$ un tore, ${\bf D}$ un groupe semi-simple et ${\bf M}$ 
une ${\bf T}\times {\bf D}$-vari{\'e}t{\'e}.
Soit ${\bf P}$ un sous-groupe parabolique de ${\bf D}$, ${\bf H}$ 
son sous-groupe d{\'e}riv{\'e} et
${\bf S}$ le centre connexe d'un sous-groupe de Levi de ${\bf P}$.

Notons $\rho\,:\,{\bf P}\times {\bf T}\longto {\bf S}\times {\bf T}$, 
l'application qui au couple 
$(p,t)\in {\bf P}\times {\bf T}$ associe le couple $(s,t)\in {\bf S}\times {\bf T}$ o{\`u} $s$ v{\'e}rifie $s^{-1}p\in {\bf H}$.

\begin{prop}
\label{prop:rho1}
Avec les notations introduites ci-dessus, on suppose que ${\bf M}$ est
factorielle et que ${\bf M}\DQD{\bf H}$ existe.

Alors, il existe un ouvert non vide $\Omega$ de ${\bf M}$ tel que 
pour tout $x$ dans $\Omega$
$$
\Ker{{\bf S}\times {\bf T}}{{\bf M}\DQD {\bf H}}
=\rho\left(
({\bf P}\times {\bf T})_x
\right).
$$
\end{prop}

\begin{preuve}
Notons $\pi\,:\,{\bf M}\longto {\bf M}\DQD {\bf H}$ l'application quotient. 
D'apr{\`e}s les lemmes~\ref{frac.inv} et~\ref{quot_ratio}, il existe un ouvert non vide
$\Omega_1$ de ${\bf M}\DQD {\bf H}$ tel que pour tout $z$ dans $\Omega_1$,
$\pi^{-1}(z)$ contient une unique orbite ouverte de ${\bf H}$.
Quitte {\`a} remplacer $\Omega_1$ par un ouvert plus petit, on peut supposer de plus que
pour tout $z$ dans $\Omega_1$, $({\bf S}\times {\bf T})_z$ est {\'e}gal {\`a}
$\Ker{{\bf S}\times {\bf T}}{{\bf M}\DQD {\bf H}}$ en vertu du 
lemme~\ref{kertore}. 

L'ensemble $\Omega_2$ des points $x$ de ${\bf M}$ tels que la dimension de
${\bf H}.x$ soit maximale est un
ouvert non vide de ${\bf M}$. Posons $\Omega=\Omega_2\cap\pi^{-1}(\Omega_1)$. 

Soit $x\in\Omega$. Posons $z=\pi(x)$.
Nous affirmons que~:
$$
({\bf S}\times {\bf T})_z=
\{
(s,t)\in {\bf S}\times {\bf T}\;:\;
(s,t).x\in {\bf H}.x
\}.
$$

En effet, $({\bf S}\times {\bf T})_z$ stabilise $\pi^{-1}(z)$ 
et permute les orbites de ${\bf H}$. 
Or, ${\bf H}.x$ est l'unique orbite ouverte de ${\bf H}$ dans $\pi^{-1}(z)$. 
Donc, $({\bf S}\times {\bf T})_z$ stabilise ${\bf H}.x$~. 

Enfin, pour tout $(s,t)\in {\bf S}\times {\bf T}$, on a~:
$$
\begin{array}{l@{\iff}l}
(s,t).x\in {\bf H}.x&\exists h\in {\bf H}\ \ \ st.x=h.x\\
&\exists h\in {\bf H}\ \ \ (hs,t).x=x\\
&(s,t)\in\rho\left(
({\bf P}\times {\bf T})_x
\right).
\end{array}
$$
\end{preuve}

\subsection{}
On conserve les notations de la section~\ref{sec:situation1}.
On suppose de plus que~:

\begin{itemize}
\item ${\bf D}$ est un sous-groupe d'un groupe semi-simple 
$\hat{{\bf D}}$.
\item ${\bf T}$ est le produit d'un tore ${\bf T}_1$ est d'un tore
  maximal $\h{\bf T}$ de $\hat{{\bf D}}$.
\item ${\bf M}$ est le produit d'une ${\bf T}_1\times {\bf D}$-vari{\'e}t{\'e} $V$ et de $\hat{{\bf D}}/\h{\bf U}$, o{\`u}
$\h{\bf U}$ est un sous-groupe unipotent maximal de $\hat{{\bf D}}$.
\item ${\bf T}_1\times\h {\bf T}\times {\bf D}$ agit sur $\hat{{\bf D}}/\h{\bf U}$ par~:
$(t,\h t,d).\h d\h{\bf U}=d\h d\h t^{-1}\h{\bf U}$.
\end{itemize}

Soit $\h {\bf B}$ le sous-groupe de Borel de $\hat{{\bf D}}$ contenant $\h{\bf U}$.
Consid{\'e}rons l'application $\hat{{\bf D}}$-{\'e}quivariante 
$q\,:\,\hat{{\bf D}}/\h{\bf U}\longto \hat{{\bf D}}/\h {\bf B}$.

Dans ce cas, la proposition~\ref{prop:rho1} est compl{\'e}t{\'e}e par la 

\begin{prop}
\label{prop:rho2}
Conservons les notations ci-dessus.
Soit $x\in X$ et $y\in \hat{{\bf D}}/\h{\bf U}$.
Alors, 
$\rho\left( ({\bf P}\times {\bf T})_{(x,y)}\right)$ est isomorphe au quotient de 
$({\bf P}\times {\bf T}_1)_{(x,q(y))}$ par son radical unipotent.
\end{prop}

\begin{preuve}
Consid{\'e}rons la projection 
$p\: :\: {\bf P}\times {\bf T}_1\times \h{\bf T}\longto {\bf P}\times {\bf T}_1$. 
Remarquons tout d'abord que comme chaque fibre de $q$ est isomorphe {\`a} $\h{\bf T}$, 
$p$ induit un isomorphisme de $({\bf P}\times {\bf T}_1\times \h{\bf
  T})_y$ sur 
$({\bf P}\times {\bf T}_1)_{q(y)}$.

Montrons par ailleurs que $\rho$ induit un isomorphisme de 
$({\bf P}\times {\bf T}_1\times \h{\bf T})_{y}/({\bf P}\times
{\bf T}_1\times \h{\bf T})_{y}^u$ sur 
$\rho(({\bf P}\times {\bf T}_1\times \h{\bf T})_{y})$.

Comme l'action de ${\bf P}$ sur $\hat{{\bf D}}/\h{\bf U}_{\bf D}^-$
est la restriction de celle de $\hat{{\bf D}}$,
${\bf P}_y$ est unipotent. Donc, ${\bf P}_y$ est inclus dans $({\bf P}\times {\bf T}_1\times \h{\bf T})_{y}^u$.
En particulier, le noyau de la restriction de $\rho$ {\`a} $({\bf P}\times {\bf T}_1\times \h{\bf T})_{y}$ est 
inclus dans $({\bf P}\times {\bf T}_1\times \h{\bf T})_{y}^u$.
L'inclusion r{\'e}ciproque est vrai car $\rho\left ( ({\bf P}\times {\bf T}_1\times \h{\bf T})_{y}^u\right)$ 
est un sous-groupe unipotent d'un tore~: il est donc trivial.

Mais alors, sur le diagramme commutatif suivant~:
\begin{diagram}
&({\bf T}_1\times\h{\bf T}\times {\bf P})_y&\rTo&({\bf T}_1\times {\bf P})_{q(y)}\\
&\dTo & &\dTo\\
{\bf T}_1\times\h{\bf T}= &\rho(({\bf T}_1\times\h{\bf T}\times {\bf P})_y)&\lTo&({\bf T}_1\times {\bf P})_{q(y)}/({\bf T}_1\times {\bf P})_{q(y)}^u,
\end{diagram}
les fl{\`e}ches horizontales sont des isomorphismes. 
On obtient alors l'isomorphisme recherch{\'e} par restriction 
{\`a} $({\bf T}_1\times\h{\bf T}\times {\bf P})_{(x,y)}$ gr{\^a}ce {\`a} l'

\noindent 
\underline{Affirmation }~: 
Soit $R$ un groupe r{\'e}soluble connexe et $R'$ un sous-groupe de $R$. 
Alors, $R'^u=R'\cap R^u$.\\

Le sous-groupe $R'\cap R^u$ est distingu{\'e} dans $R'$ et unipotent~:
il est donc inclus dans $R'^u$. 

Par ailleurs, $R'/(R'\cap R^u)$ s'injecte dans $R/R^u$ qui est un tore. 
Donc, $R'/(R'\cap R^u)$ est r{\'e}ductif et $ (R'\cap R^u)\subset R'^u$.
\end{preuve}

\bibliographystyle{smfalpha}
\bibliography{biblio}

\begin{center}
  -\hspace{1em}$\diamondsuit$\hspace{1em}-
\end{center}

\end{document}